\documentclass{amsart}
\usepackage{amscd,amsthm,amssymb,amsfonts,amsmath,euscript}

\theoremstyle{plain}
\newtheorem{thm}{Theorem}[section]
\newtheorem{lemma}[thm]{Lemma}
\newtheorem{prop}[thm]{Proposition}

\theoremstyle{definition}

\theoremstyle{remark}
\newtheorem{remark}[thm]{Remark}

\newcommand{\nc}{\newcommand}

\def\makeop#1{\expandafter\def\csname#1\endcsname
  {\mathop{\rm #1}\nolimits}\ignorespaces}
\makeop{Hom}   \makeop{End}   \makeop{Aut}   \makeop{Isom}  \makeop{Pic} 
\makeop{Gal}   \makeop{ord}   \makeop{Char}  \makeop{Div}   \makeop{Lie} 
\makeop{PGL}   \makeop{Corr}  \makeop{PSL}   \makeop{sgn}   \makeop{Spf}
\makeop{Spec}  \makeop{Tr}    \makeop{Nr}    \makeop{Fr}    \makeop{disc}
\makeop{Proj}  \makeop{supp}  \makeop{ker}   \makeop{im}    \makeop{dom}
\makeop{coker} \makeop{Stab}  \makeop{SO}    \makeop{SL}    \makeop{SL}
\makeop{Cl}    \makeop{cond}  \makeop{Br}    \makeop{inv}   \makeop{rank}
\makeop{id}    \makeop{Fil}   \makeop{Frac}  \makeop{GL}    \makeop{SU}
\makeop{Trd}   \makeop{Sp}    \makeop{Tr}    \makeop{Trd}   \makeop{diag}
\makeop{Res}   \makeop{ind}   \makeop{depth} \makeop{Tr}    \makeop{st}
\makeop{Ad}    \makeop{Int}   \makeop{tr}    \makeop{Sym}   \makeop{can}
\makeop{length}\makeop{SO}    \makeop{torsion} \makeop{GSp} \makeop{Ker}
\def\makebb#1{\expandafter\def
  \csname bb#1\endcsname{{\mathbb{#1}}}\ignorespaces}
\def\makebf#1{\expandafter\def\csname bf#1\endcsname{{\bf
      #1}}\ignorespaces} 
\def\makegr#1{\expandafter\def
  \csname gr#1\endcsname{{\mathfrak{#1}}}\ignorespaces}
\def\makescr#1{\expandafter\def
  \csname scr#1\endcsname{{\EuScript{#1}}}\ignorespaces}
\def\makecal#1{\expandafter\def\csname cal#1\endcsname{{\mathcal
      #1}}\ignorespaces} 

\def\doLetters#1{#1A #1B #1C #1D #1E #1F #1G #1H #1I #1J #1K #1L #1M
                 #1N #1O #1P #1Q #1R #1S #1T #1U #1V #1W #1X #1Y #1Z}
\def\doletters#1{#1a #1b #1c #1d #1e #1f #1g #1h #1i #1j #1k #1l #1m
                 #1n #1o #1p #1q #1r #1s #1t #1u #1v #1w #1x #1y #1z}
\doLetters\makebb   \doLetters\makecal  \doLetters\makebf
\doLetters\makescr 
\doletters\makebf   \doLetters\makegr   \doletters\makegr
     \def\qed{\qedmark\medbreak}%
\def\qedmark{{\enspace\vrule height 6pt width 5pt depth 1.5pt}}%

\normalsize

\makeop{Bl}

\def\Spec{{\rm Spec}}

\def\Fpbar{\overline{\bbF}_p}
\def\Fp{{\bbF}_p}
\def\Fq{{\bbF}_q}

\def\Qpbar{\overline{{\bbQ}_p}}
\def\Qp{{\bbQ}_p}

\def\Zp{{\bbZ}_p}
\def\Qbar{\overline{\bbQ}}

\newcommand{\Z}{\mathbb Z}
\newcommand{\Q}{\mathbb Q}
\newcommand{\R}{\mathbb R}
\newcommand{\C}{\mathbb C}
\renewcommand{\H}{\mathbb H}  
\newcommand{\A}{\mathbb A}    
\newcommand{\F}{\mathbb F}

\newcommand{\npr}{\noindent }

\newcommand{\<}{\langle}   
\renewcommand{\>}{\rangle} 

\nc{\embed}{\hookrightarrow}

\newcommand{\ch}{characteristic }
\newcommand{\ac}{algebraically closed }
\newcommand{\dieu}{Dieudonn\'{e} }

\nc{\ol}{\overline}
\nc{\wt}{\widetilde}
\nc{\opp}{\mathrm{opp}}
\def\ul{\underline}

\makeop{Ram}
\makeop{Rep}

\begin{document}
\renewcommand{\thefootnote}{\fnsymbol{footnote}}
\setcounter{footnote}{-1}
\numberwithin{equation}{section}


\title{An exact geometric mass formula} 
\author{Chia-Fu Yu}
\address{
Institute of Mathematics \\
Academia Sinica \\
128 Academia Rd.~Sec.~2, Nankang\\ 
Taipei, Taiwan \\ and NCTS (Taipei Office)}
\email{chiafu@math.sinica.edu.tw}
\address{
Max-Planck-Institut f\"ur Mathematik \\
Vivatsgasse 7 \\
Bonn, 53111\\ 
Germany}

\date{June 19, 2007.  The research is partially supported by NSC
 96-2115-M-001-001.}
\begin{abstract}
We show an exact geometric mass formula for superspecial points in the
reduction of any quaternionic Shimura variety modulo at a good prime $p$.
\end{abstract} 

\maketitle


\section{Introduction}
\label{sec:01}
\def\Mass{{\rm Mass}}

Let $p$ be a rational prime number. 
Let $B$ be a totally indefinite quaternion algebra 
over a totally real field $F$ of degree $d$, together with a positive
involution $*$. Assume that $p$ is unramified
in $B$. Let $O_B$ be a maximal order stable under the
involution $*$. Let $(V,\psi)$ be a non-degenerate $\Q$-valued
skew-Hermitian (left) $B$-module with dimension $2g$ over
$\Q$. Put $m:=\frac{g}{2d}$, a positive integer. 
A polarized abelian $O_B$-variety $\ul A=(A,\lambda,\iota)$ is a polarized
abelian variety $(A,\lambda)$ together with a ring 
monomorphism $\iota:O_B\to \End(A)$ such that 
$\lambda\circ \iota(b^*)=\iota(b)^t\circ \lambda$ for all
$b\in O_B$. Let $k$ be an \ac field of \ch $p$. An abelian variety
over $k$ is said to be superspecial if it is isomorphic to a product of
supersingular elliptic curves. Denote by $\Lambda^B_{g}$ the set of
isomorphism classes of $g$-dimensional superspecial principally
polarized abelian $O_B$-varieties over $k$. Define the mass of
$\Lambda^{B}_{g}$ to be 
\begin{equation}
  \label{eq:11}
  \Mass(\Lambda^{B}_{g}):=\sum_{\ul A\in \Lambda^B_{g}}
  \frac{1}{|\Aut(A,\lambda,\iota)|}. 
\end{equation}

The mass ${\rm Mass}(\Lambda^{B}_{g})$ is studied in Ekedahl~\cite{ekedahl:ss}
(Ekedahl's result relies on an explicit volume computation in
Hashimoto-Ibukiyama~\cite[Proposition 9, 
p.~568]{hashimoto-ibukiyama:classnumber}) in the
special case $B=M_2(\Q)$. He proved

\begin{thm} [Ekedahl, Hashimoto-Ibukiyama]\label{11} One has
\begin{equation}
  \label{eq:12}
  \Mass(\Lambda_{g})=
  \frac{(-1)^{g(g+1)/2}}{2^g}  \prod_{i=1}^g \zeta(1-2i)
  \cdot \prod_{i=1}^{g} p^i+(-1)^i, 
\end{equation}
where $\Lambda_{g}$ is the set of isomorphism classes of
$g$-dimensional superspecial principally polarized abelian varieties
over $k$ and $\zeta(s)$ is the Riemann zeta function.
\end{thm}

Let $B_{p,\infty}$ be the quaternion algebra over $\Q$ ramified
exactly at $\{p,\infty\}$. Let $B'$ be the
quaternion algebra over $F$ such that
$\inv_v(B')=\inv_v(B_{p,\infty}\otimes_\Q B)$ for all $v$. Let
$\Delta'$ be the discriminant of $B'$ over $F$. 

In this paper we prove
\begin{thm}\label{12} One has
  \begin{equation}\label{eq:13}
\Mass(\Lambda^B_g)=\frac{(-1)^{dm(m+1)/2}}{2^{md}}
\prod_{i=1}^m \left \{ 
  \zeta_F(1-2i)   
\prod_{v|\Delta'}
N(v)^i+(-1)^i \prod_{v|p, v\nmid\Delta'} N(v)^i+1 \right \},
\end{equation}
where $\zeta_F(s)$ is the Dedekind zeta function. 
\end{thm}

Let $N\ge 3$ be a prime-to-$p$ positive integer. Choose a primitive
$n$-th root of unity $\zeta_N\in \Qbar\subset \C$ and fix an
embedding $\Qbar \hookrightarrow \Qpbar$. Let $\calM$ be the moduli
space over $\Fpbar$ of $g$-dimensional principally polarized abelian
$O_B$-varieties with a symplectic $O_B$-linear level-$N$ structure
w.r.t. $\zeta_N$. Let $L_0$ be a self-dual $O_B$-lattice of $V$ with
respect to $\psi$. 
Let $G_1$ be the automorphism group scheme over $\Z$
associated to the pair $(L_0,\psi)$. As an immediate consequence of 
Theorem~\ref{12}, we get

\begin{thm}\label{13}
  The moduli space $\calM$ has
  \begin{equation}
    \label{eq:14}
    |G_1(\Z/N\Z)|\frac{(-1)^{dm(m+1)/2}}{2^{md}}
\prod_{i=1}^m \left \{ 
  \zeta_F(1-2i)   
\prod_{v|\Delta'}
N(v)^i+(-1)^i \prod_{v|p, v\nmid\Delta'} N(v)^i+1 \right \}
  \end{equation}
superspecial points. 
\end{thm}

We divide the proof of Theorem~\ref{12} into 4 parts; each part is
treated in one section. The first part is to express
the weighted sum in terms of an arithmetic mass; this is done in the
author's recent work \cite{yu:smf}. The second part is to compute the 
mass associated to a quaternion unitary group and a standard open
compact subgroup; this is done by Shimura \cite{shimura:1999}
(re-obtained by Gan and 
J.-K.~Yu \cite[11.2, p.~522]{gan-yu}) using the theory
of Bruhat-Tits Buildings). 
The third part is to compare the derived arithmetic mass in
Section 1 with ``the'' standard mass in Section 2. 
This reduces the problem to computing
a local index at $p$. The last part uses \dieu theory to compute this
local index. A crucial step is choosing a good basis for the
superspecial \dieu module concerned; this makes the computation easier.\\



\def\DS{{\rm DS}}

\npr{\bf Notation.} $\bbH$ denotes the Hamilton quaternion algebra
over $\R$. $\A_f$ denotes the finite adele ring of $\Q$ and $\hat
\Z=\prod_p \Zp$. For a number field $F$ and a finite place $v$, denote
by $O_F$ the ring of integers, $F_v$ the completion of $F$ at $v$,
$e_v$ the ramification index for $F/\Q$, $\kappa_v$ the 
residue field, $f_v:=[\kappa_v:\Fp]$ and $q_v:=N(v)=|\kappa_v|$. 
For an $O_F$-module $A$, write $A_v$ for
$A\otimes_{O_F} {O_{F,v}}$. For a scheme $X$ over $\Spec\, A$ and an
$A$-algebra $B$, write $X_B$ for $X\times_{\Spec A} \Spec B$.  
For a linear algebraic group $G$ over
  $\Q$ and an open compact subgroup $U$ of $G(\A_f)$, denote by $\DS
  (G,U)$ the double coset space $G(\Q)\backslash G(\A_f)/U$, and write
  $\Mass(G,U):=\sum_{i=1}^h |\Gamma_i|^{-1}$ if $G$ is
  $\R$-anisotropic, where $\Gamma_i:=G(\Q)\cap c_i U c_i^{-1}$ and
  $c_1,\dots,c_h$ are complete representatives for $\DS(G,U)$. For a
  central simple algebra $B$ over $F$, write $\Delta(B/F)$ for the
  discriminant of $B$ over $F$.
If $B$ a central division algebra over a non-archimedean 
local field $F_v$, denote by $O_B$ the maximal order of $B$, $\grm(B)$ the
maximal ideal and $\kappa(B)$ the residue field. $\Q_{p^n}$ denotes
the unramified extension of $\Q_p$ of degree $n$ and 
write $\Z_{p^n}:=O_{\Q_{p^n}}$. 


\section{Simple mass formulas}
\label{sec:02}
\def\qisom{{\rm Q\text{-}isom}}
\newcommand{\double}[2]{#1(\Q)\backslash #1(\A_f)/#2}
\newcommand{\doublecoset}[1]{\double{#1}{#1(\hat \Z)}}

Let $B$ be a finite-dimensional semi-simple algebra over $\Q$ with a
positive involution $*$, and $O_B$ be an order of $B$ stable under
$*$. Let $k$ be any field.

To any polarized abelian $O_B$-varieties $\ul A=(A,\lambda,\iota)$
over $k$,
we associate a pair $(G_x,U_x)$, where $G_x$ is the group scheme over
$\Z$ representing the functor 
\[ R \mapsto \{h\in (\End_{O_B}(A_k)\otimes R)^\times\, |\, h' h=1\},
\]
where $h\mapsto h'$ is the Rasoti involution, and $U_x$ is the open
compact subgroup $G_x(\hat \Z)$.
For any prime $\ell$, we write $\ul A(\ell)$ for the associated
$\ell$-divisible group with additional structures $(A[\ell^\infty],
\lambda_\ell, \iota_\ell)$, where $\lambda_\ell$ is the induced
quasi-polarization from $A[\ell^\infty]$ to
$A^t[\ell^\infty]=A[\ell^\infty]^t$ (the Serre dual), 
and $\iota_\ell:O_B\otimes\Z_\ell\to \End(A[\ell^\infty])$ 
the induced ring monomorphism. 
For any two objects $\ul A_1$ and $\ul A_2$ over $k$,
denote by $\qisom_k(\ul A_1, \ul A_2)$  the set of $O_B$-linear 
quasi-isogenies $\varphi:A_1\to A_2$ over $k$ such that
$\varphi^*\lambda_2=\lambda_1$, and   
$\Isom_k(\ul A_1(\ell), \ul A_2(\ell))$ the set of 
$O_B\otimes \Z_\ell$-linear isomorphisms $\varphi:A_1[\ell^\infty]\to
A_2[\ell^\infty]$ over $k$ such that $\varphi^*\lambda_2=\lambda_1$.

Let $x:=\ul A_0=(A_0,\lambda_0,\iota_0)$ be a fixed 
polarized abelian $O_B$-variety over $k$.  
Denote by $\Lambda_x(k)$ the set of isomorphisms classes of polarized
abelian $O_B$-varieties $\ul A$ over $k$ such that 
\begin{itemize}
\item [($I_\ell$):] $\Isom_k(\ul
A_0(\ell),\ul A(\ell))\neq \emptyset$ for all primes $\ell$. 
\end{itemize}
Let $\Lambda'_x(k)\subset \Lambda_x(k)$ be the subset consisting of
objects such that 
\begin{itemize}
\item [(Q):] $\qisom_k(\ul A_0,\ul A)\neq \emptyset$.
\end{itemize}
Let $\ker^1(\Q,G_x)$ denote the kernel of the local-global map
  $H^1(\Q,G_x)\to \prod_v  H^1(\Q_v,G_x)$.

\begin{thm} {\bf (\cite[Theorem~2.3]{yu:smf})}
\label{21}  Suppose that $k$ is
  a field of finite type over 
  its prime field.  

\npr {\rm (1)} There is a natural bijection  
$\Lambda'_x(k)\simeq  \DS(G_x,U_x).$
Consequently, $\Lambda'_x(k)$ is finite.

\npr {\rm (2)} One has $\Mass(\Lambda'_x(k))=\Mass(G_x,U_x)$.
\end{thm}


\begin{thm}\label{22} {\bf (\cite[Theorem~4.6 and Remark~4.7]{yu:smf})} 
  Notation as above. 
  If $k\supset \Fp$ is \ac and $A_0$ is supersingular, then
  $\Mass(\Lambda_x'(k))=\Mass(G_x,U_x)$ and
  $\Mass(\Lambda_x(k))=|\ker^1(\Q,G_x)|\cdot\Mass(G_x,U_x)$.  
\end{thm}

\begin{remark}
   The statement of Theorem~\ref{22} 
is valid for basic abelian $O_B$-varieties
in the sense of Kottwitz (see \cite{rapoport-zink} for the definition). 
The present form is enough for our purpose. 
\end{remark}

\section{An exact mass formula of Shimura}
\label{sec:03}

Let $D$ be a totally definite quaternion division algebra over a
totally real field $F$ of degree $d$. Let (bar) $d\mapsto \bar d$
denote the canonical involution. 
Let $(V',\varphi)$ be a $D$-valued totally definite 
quaternion Hermitian $D$-module of rank $m$. 
Let $G^\varphi$ 
denote the unitary group attached to $\varphi$. This is a reductive
group over $F$ and is regarded as a group over $\Q$ via the
Weil restriction of scalars from $F$ to $\Q$. 
Choose a maximal order $O_D$ of $D$ stable under the canonical
involution $ \bar{\ } $. 
Let $L$ be an $O_D$-lattice in $V'$ which is maximal among the lattices 
on which $\varphi$ takes its values in $O_B$. Let $U_0$ be the 
open compact subgroup of $G^\varphi(\A_f)$ which stabilizes the adelic
lattice $L\otimes_\Z \hat{\Z}$.

The following is deduced from a mass formula of
Shimura \cite{shimura:1999} (also see 
Gan - J.-K.~Yu \cite[11.2, p.~522]{gan-yu}). This form is more
applicable to prove Theorem~\ref{12}.   

\begin{thm} [Shimura] \label{31} One has
  \begin{equation}\label{eq:31}
\Mass(G^\varphi, U_0)=\frac{(-1)^{dm(m+1)/2}}{2^{md}}
\prod_{i=1}^m \left \{ 
  \zeta_F(1-2i)   
\prod_{v|\Delta(D/F)}
N(v)^i+(-1)^i\right \}.
\end{equation}
\end{thm}

\npr Deduction. In \cite[Introduction, p.~68]{shimura:1999} 
Shimura gives the explicit formula
\begin{equation}\label{eq:32}
 \Mass(G^\varphi, U_0)=|D_F|^{m^2}  \prod_{i=1}^m D_F^{1/2}\left 
[(2i-1)!(2\pi)^{-2i}\right ]^d \zeta_F(2i) \cdot \prod_{v|\Delta(D/F)}
\prod_{i=1}^{m} N(v)^i+(-1)^i,
\end{equation}
where $D_F$ is the discriminant of $F$ over $\Q$. 
Using the functional equation for $\zeta_F(s)$, we deduce
(\ref{eq:31}) from (\ref{eq:32}).



\section{Global comparison}
\label{sec:04}

Keep the notation as in Section~\ref{sec:01}. Fix a $g$-dimensional
superspecial principally polarized abelian $O_B$-variety
$x=(A_0,\lambda_0,\iota_0)$ over $k$. Define $\Lambda_x:=\Lambda_x(k)$
as in Section~\ref{sec:02}. Let $(G_x,U_x)$ be the pair associated to
$x$.

\begin{lemma}\label{41}
  Any two self-dual $O_B\otimes \Zp$-lattices of $(V_{\Qp},\psi)$ are
  isomorphic. 
\end{lemma}
\begin{proof}
  The proof is elementary and omitted.
\end{proof}

\begin{lemma}\label{42}
  One has (1) $\Lambda_x=\Lambda_g^B$ (2) $\ker^1(\Q,G_x)=\{1\}$.
\end{lemma}
\begin{proof}
  (1) The inclusion $\Lambda_x\subset \Lambda_g^B$ is clear. We show
      the other direction. Let $\ul A\in \Lambda^B_g$. It follows from
      Lemma~\ref{41} that the condition $(I_\ell)$ is satisfied for
      primes $\ell\neq p$. Let $M$ be the covariant \dieu module of
      $A$. One chooses an isomorphism $O_{B,p}\simeq M_2(O_{F,p})$ so
      that $*:(a_{ij})\mapsto (a_{ij})^t$. Using the Morita
      equivalence, it suffices to show that any two superspecial
      principally quasi-polarized \dieu modules with compatible
      $O_{F,p}$-action 
      are isomorphic. This follows from Theorem~\ref{51}.

  (2) Since $G_x$ is semi-simple and simply connected (as it is an
      inner form of $\Res_{F/\Q} \Sp_{2m,F}$), the Hasse principle for $G_x$
      holds.  \qed 
\end{proof}

\subsection{}
\label{sec:41}
We compute that
\begin{itemize}
\item [(i)] $G_x(\R)=\{h\in M_m(\H)^d \, | \, \bar h^t h=1\,\}$,
\item [(ii)] for $\ell\neq p$, we have $G_x(\Q_\ell)=\prod_{v|\ell} G_{x,v}$
  and $U_{x,\ell}=\prod_{v|\ell} U_{x,v}$,
  where
\begin{equation}
  \label{eq:41}
  \begin{split}
  &  G_{x,v}=
  \begin{cases}
  \Sp_{2m}(F_v), & \text{if $v\nmid\Delta(B/F)$}, \\
  \{h\in M_m(B_v)\, | \, \bar h^t h=1\,\}, & \text{otherwise,}
  \end{cases} \\
  & U_{x,v}=
  \begin{cases}
  \Sp_{2m}(O_{F_v}), & \text{if $v\nmid\Delta(B/F)$}, \\
  \{h\in M_m(O_{B_v})\, | \, \bar h^t h=1\,\}, & \text{otherwise},
  \end{cases}
  \end{split}
\end{equation}
\item [(iii)]  $G_x(\Q_p)=\prod_{v|p} G_{x,v}$, where 
  \begin{equation}
  \label{eq:42}
   G_{x,v}=
  \begin{cases}
  \Sp_{2m}(F_v), & \text{if $v\nmid\Delta'$}, \\
  \{h\in M_m(B'_v)\, | \, \bar h^t h=1\,\}, & \text{otherwise.}
  \end{cases}.
\end{equation}
\end{itemize}

Take $D=B'$ and $V'=D^{\oplus m}$ with $\varphi(\ul x,\ul y)=\sum x_i \bar
  y_i$, and take $L=O_D^{\oplus m}$. We compute that
\begin{itemize}
\item [(i)'] $G^\varphi(\R)=\{h\in M_m(\H)^d \, | \, \bar h^t h=1\,\}$,
\item [(ii)'] for any $\ell$, we have 
$G_x(\Q_\ell)=\prod_{v|\ell} G^\varphi_{v}$
  and $U_{0,\ell}=\prod_{v|\ell} U_{0,v}$,
  where
\begin{equation}
  \label{eq:43}
  \begin{split}
  &  G^\varphi_{v}=
  \begin{cases}
  \Sp_{2m}(F_v), & \text{if $v\nmid\Delta'$}, \\
  \{h\in M_m(B'_v)\, | \, \bar h^t h=1\,\}, & \text{otherwise,}
  \end{cases} \\
  & U_{0,v}=
  \begin{cases}
  \Sp_{2m}(O_{F_v}), & \text{if $v\nmid\Delta'$}, \\
  \{h\in M_m(O_{B'_v})\, | \, \bar h^t h=1\,\}, & \text{otherwise.}
  \end{cases}
  \end{split}
\end{equation}
\end{itemize}
For $\ell\neq p$ and $v|\ell$, one has $B_v=B_v'$ and that
$v\nmid\Delta(B/F)$ if and only if $v\nmid\Delta'$.
It follows from computation above that $G_{x,\R}\simeq G^\varphi_\R$ and 
$G_{x,{\Q_\ell}}\simeq G^\varphi_{\Q_\ell}$ for all $\ell$. Since the
Hasse principle holds for the adjoint group $G_x^{\rm ad}$, we get
$G_x\simeq G^\varphi$ over $\Q$. We fix an isomorphism and write
$G_x=G^\varphi$. For $\ell\neq p$ and $v|\ell$, the subgroups
$U_{0,v}$ and $U_{x,v}$ are conjugate, and hence they have the same
local volume.

\subsection{}
\label{sec:42}
Applying Theorem~\ref{22} in our setting (Section 1) and using
Lemma~\ref{42}, we get $\Mass(\Lambda_g^B)=\Mass(G_x,U_x)$. Using 
the result in Subsection~\ref{sec:41}, we get
\begin{equation}
  \label{eq:44}
  \Mass(\Lambda_g^B)=\Mass(G^\varphi,U_0)\cdot \mu(U_{0,p}/U_{x,p}),
\end{equation}
where $\mu(U_{0,p}/U_{x,p})=[U_{x,p}:U_{0,p}\cap
U_{x,p}]^{-1}[U_{0,p}:U_{0,p}\cap U_{x,p}]$.   

\section{Local index $\mu(U_{0,p}/U_{x,p})$}
\label{sec:05}

Let $(M',\<\,,\>',\iota')$ be the covariant \dieu module associated to
the point $x=(A_0,\lambda_0,\iota_0)$ in the previous section. Choose
an isomorphism $O_B\otimes \Z_p\simeq M_2(O_F\otimes \Z_p)$ so that
$*$ becomes the transpose. Let $M:=eM'$, $\<\, ,\>:=\<\, ,\>'|_{M}$
and $\iota:=\iota'|_{O_F}$, where $e=
\begin{pmatrix}
  1 & 0 \\ 0 & 0
\end{pmatrix}$ in $M_2(O_F\otimes \Z_p)$. The triple $(M, \<\, ,\>,\iota)$ is a
superspecial principally quasi-polarized \dieu module with compatible
$O_F\otimes \Z_p$-action of rank $g=2dm$. Let $M=\oplus_{v|p} M_v$ be the
decomposition with respect to the decomposition $O_F\otimes
\Z_p=\oplus_{v|p} \calO_v$; here we write $\calO_v$ for $O_{F_v}$. By
the Morita equivalence, we have
\begin{equation}
  \label{eq:51}
  U_{x,p}=\Aut_{{\rm DM},O_B}(M',\<\, ,\>')=
  \Aut_{{\rm DM},O_F}(M,\<\, ,\>)=\prod_{v|p} U_{x,v},
\end{equation}
where $U_{x,v}:= \Aut_{{\rm DM},\calO_v}(M_v,\<\, ,\>)$. 

\def\Zf{{\Z/f_v\Z}}

Let $W:=W(k)$ be ring of Witt vectors over $k$ and $\sigma$ the
absolute Frobenius map on $W$. Let $\scrI:=\Hom(\calO_v,W)$ be the set of
embeddings; write $\scrI=\{\sigma_i\}_{i\in \Z/f_v\Z}$ so that $\sigma
\sigma_i=\sigma_{i+1}$ for all $i$. We identify $\Z/f_v\Z$ with $\scrI$ through
$i\mapsto \sigma_i$. Decompose $M_v=\oplus_{i\in \Zf} M_v^i$ into
$\sigma_i$-isotypic components $M^i_v$. One has (1) each component
$M^i_v$ is a free $W$-module of rank $2m$, which is self-dual with respect to the
pairing $\<\, ,\>$, (2) $\<M^i_v,M^j_v\>=0$ if $i\neq j$, and 
(3) the operations $F$ and $V$ shift by degree 1 and degree -1, respectively.

\begin{thm}\label{51}
  Let $(M_v, \<\, ,\>,\iota)$ be as above. There is a symplectic basis
  $\{X^i_j,Y^i_j\}_{j=1,\dots, m}$ for $M^i_v$ such that 
  \begin{itemize}
  \item [(i)]  $Y^i_j\in V M_v^{i+1}$,
  \item [(ii)] $FX^i_j=-Y^{i+1}_j$ and $FY^i_{j}=pX^{i+1}_j$,
  \end{itemize}
 for all $i\in \Zf$ and all $j$. 
\end{thm}
\begin{proof}
We write $f$, $M$ and $q$ for $f_v$, $M_v$ and $q_v$, respectively. 
Suppose that $f=2c$ is even. Let $N:=\{x\in M \, |\, F^c
  x=(-1)^c V^c x \}$. Since $M$ is superspecial, we have ($*$) $F^2
  {N}=p {N}$, $\wt{N}\otimes_{\Z_{q}} W\simeq M$ and $N=\oplus N^i$. 
Since $\ol{VN^1}$ is isotropic
with respect to $\<\, ,\>$ in $N/pN$, we can choose a symplectic
basis $\{X^0_j,Y^0_j\}_{j=1,\dots,m}$ for $N^0$ such that $Y^0_j\in
VN^1$ for all $j$.
Define $X^i_j$ and $Y^i_j$ recursively for $j=1,\cdots, j$:
\begin{equation}\label{rec}
  X^i_{j+1}=p^{-1}FY^i_j,\quad Y^i_{j+1}=-F X^i_j.
\end{equation}
One has $X^i_{j+2}=\frac{-1}{p}F^2 X_j^i$ and $\ Y^i_{j+2}=\frac{-1}{p} F^2
  Y_j^i$; hence
\[  X^f_j=(-1)^c p^{-c} F^{2c} X^0_j=X^0_j, \quad
 Y^f_j=(-1)^c p^{-c} F^{2c} Y^0_j=Y^0_j, \]
for all $j$. It is easy to see that $\{X^i_j,Y^i_j\}_{j=1,\dots, m}$
forms a symplectic basis for $N^i$. 

Suppose that $f=2c+1$ is odd. Let $N:=\{
  x\in M\, |\, F^{2f}x+p^f x=0\}$. We construct a symplectic basis
  $\{X^0_j,Y^0_j\}_{j=1,\dots, m}$ for $N^0$ with the properties:
  $X^0_j\not\in VN^1$, $Y^0_j\in VN^1$ and $Y^0_j=(-1)^{c+1}p^{-c}F^f
  X^0_j$ for all $j$. 
We can choose $X^0_1\in N^0 \backslash VN^1$ so that
  $\<X^0_1,(-1)^{c+1}p^{-c}F^f X^0_1\>\in \Z_{q^2}^\times$. This
  follows from the fact that the form $(x,y):=\<x,p^{-c}F^f y\> \mod
  p$ is a non-degenerate Hermitian form on $N^0/VN^1$. 
Set $Y^0_1=(-1)^{c+1}p^{-c}F^f X^0_1$ and let
  $\mu:=\<X^0_1,Y^0_1\>$. From $\<F^f X^0_1,.F^f
  Y^0_1\>=\<(-1)^{c+1}p^c Y^0_1,(-1)^c p^{c+1} X^0_1\>$, we get
  $\mu\in \Z^\times_q$. Since $\Q_{q^2}/\Q_{q}$ is unramified,
  replacing $X^0_1$ by a suitable $\lambda X^0_1$, we get
  $\<X^0_1,Y^0_1\>=1$. Do the same construction for the complement of
  the submodule $<X^0_1, Y^0_1>$ and use induction; we exhibit such a
  basis for $N^0$.

Define $X^i_j$ and $Y^i_j$ recursively for $i=1,\dots, f$ as
(\ref{rec}). We verify again that $X^f_j=X^0_j$ and $Y^f_j=Y^0_j$. 
It follows from the relation (\ref{rec}) that 
$\{X^i_j,Y^i_j\}_{j=1,\dots, m}$ forms a
symplectic basis for $N^i$ for all $i$. This completes the
proof. \qed 
\end{proof}

\begin{prop} \label{52} Notation as above.

{\rm (1)} If $f_v$ is even, then 
\begin{equation}\label{eq:53}
  U_{x,v}=\left \{
  \begin{pmatrix}
    A & B \\ C & D
  \end{pmatrix}\in \Sp_{2m}(\Z_{q_v})\, |\, B \equiv 0 \mod p\, \
\right \}. 
\end{equation}

{\rm (2)} If $f_v$ is odd, then 
\begin{equation}
  \label{eq:54}
   U_{x,v}\simeq \{h\in M_m(O_{B'_v})\,
| \, \bar h^t h=1\,\}.
\end{equation}
\end{prop}

\begin{proof}
Let $\phi\in U_{x,v}$. Choose a symplectic basis $\calB$ for $M_v$ as in
Theorem~\ref{51}. Since $\phi$ commutes with the $O_F$-action, we have
$\phi=(\phi_i)$, where $\phi_i\in \Aut(M^i_v,\<\,,\>)$. Write
$\phi_i=\begin{pmatrix} 
    A_i & B_i \\ C_i & D_i
\end{pmatrix}\in \Sp_{2m}(W)$ using the basis $\calB$. Since the map
$F$ is 
injective, $\phi_0$ determines the remaining $\phi_i$.  From $\phi
F^2=F^2 \phi$, we have $\phi_{i+2}=\phi_i^{(2)}$ (as matrices). Here
we write 
$\phi_i^{(n)}$ for $\phi_i^{\sigma^n}$. From $\phi F=F \phi$ we get
$A_i^{(1)}=D_{i+1}$, $B_i^{(1)}=-pC_{i+1}$, $pC_i^{(1)}=-B_{i+1}$ and
$D_i^{(1)}=A_{i+1}$.

(1) If $f_v$ is even, then $A_0,B_0,C_0,D_0\in \Z_{q_v}$ and $B_0\equiv
    0\mod p$. This shows (\ref{eq:53}).

(2) Suppose $f_v$ is odd. From $\phi_0^{(f_v+1)}=\phi_1$ we get
    $A_0^{(f_v)}=D_0$, $B_0^{(f_v)}=-pC_0$, $pC_0^{(f_v)}=-B_0$,
    $D_0^{(f_v)}=A_0$. Hence
\[ U_{x,v}=\left\{ 
   \begin{pmatrix} 
    A & -pC^\tau \\ C & A^\tau
\end{pmatrix}\in \Sp_{2m}(\Z_{q_v^2})\right\}, \]
where $\tau$ is the involution of $\Q_{q_v^2}$ over $\Q_{q_v}$. Note
that $O_{B_v'}=\Z_{q_v^2}[\Pi]$ with $\Pi^2=-p$ and $\Pi a=a^\tau
\Pi$ for all $a\in \Z_{q_v^2}$. The map $A+C\Pi \mapsto \begin{pmatrix} 
    A & -pC^\tau \\ C & A^\tau
\end{pmatrix}$ gives rise to an isomorphism (\ref{eq:54}). This proves
the proposition. \qed
\end{proof}
Let $(V_0=\F_q^{2m},\psi_0)$ be a standard symplectic space. 
Let $P$ be the stabilizer of the standard maximal isotropic subspace
$\Fq <e_1,\dots,e_m>$.
\begin{lemma}\label{53}
$ |\Sp_{2m}(\Fq)/P|= \prod_{i=1}^m (q^i+1).$
\end{lemma}
\begin{proof}
We have a natural bijection between the
group $\Sp_{2m}(\Fq)$ and the set $\calB(m)$ of
ordered symplectic bases $\{v_1,\dots,v_{2m}\}$ for $V_0$.
The first vector $v_1$ has $q^{2m}-1$ choices. The first companion 
vector $v_{m+1}$ has $q^{2m-1}$ choices as it 
does not lie in the hyperplane $v_1^{\bot}$ and we require 
$\psi_0(v_1,v_{m+1})=1$.
The remaining ordered symplectic basis can be chosen from the
complement $\Fq<v_1,v_{m+1}>^{\bot}$. Therefore, we have proved the
recursive formula 
$ |\Sp_{2m}(\Fq)|=(q^{2m}-1)q^{2m-1}|\Sp_{2m-2}(\Fq)|$. 
From this, we get
\begin{equation}
  \label{eq:55}
  |\Sp_{2m}(\Fq)|=q^{m^2}\prod_{i=1}^m (q^{2i}-1).
\end{equation}
We have
\[ P=\left\{
    \begin{pmatrix}
      A & B \\ 0 & D \\
    \end{pmatrix} ; AD^t=I_m,\ BA^t=AB^t \right\}. \]
This yields
\begin{equation}
  \label{eq:56}
  |P|=q^{\frac{m^2+m}{2}}\, |\GL_m(\Fq)|=q^{m^2} 
\prod_{i=1}^m (q^i-1).
\end{equation}
From (\ref{eq:55}) and (\ref{eq:56}), we prove the lemma.\qed
\end{proof}

By Proposition~\ref{52} and Lemma~\ref{53}, we get
\begin{thm}\label{54} One has
  \begin{equation}
\label{eq:57}
\mu(U_{0,p}/U_{x,p})=\prod_{v|p}\mu(U_{0,v}/U_{x,v})=\prod_{v|p,
  v\nmid \Delta'} \prod_{i=1}^m (q_v^i+1).    
\end{equation}
\end{thm}

Plugging the formula (\ref{eq:57}) in the formula (\ref{eq:44}), we
get the formula (\ref{eq:13}). The proof of Theorem~\ref{12} is
complete. \\

\npr{\it Acknowledgments.}
The present work relies on Shimura's paper \cite{shimura:1999} and
is also inspired by  W.-T.~Gan and J.-K.~Yu's paper \cite{gan-yu}. It
is a great pleasure to thank them. 


\end{document}